\documentclass[aop,citesort,MSNbibl,dvips]{arximspdf}

\doi{10.1214/10-AOP631}
\volume{40}
\issue{2}
\pubyear{2012}
\firstpage{743}
\lastpage{764}

\makeatletter

\newcommand{\calL}{\mathcal{L}}

\newcommand{\hT}{\hat{T}}

\newcommand{\tilf}{\tilde{f}}
\newcommand{\ra}{\rightarrow}
\newcommand{\Z}{\mathbb Z}
\newcommand{\E}{\mathbb E}
\newcommand{\R}{\mathbb R}
\newcommand{\N}{\mathbb N}
\newcommand{\PP}{\mathbb P}
\newcommand{\al}{\alpha}
\newcommand{\et}{\eta}
\newcommand{\be}{\beta}
\newcommand{\ep}{\varepsilon}
\newcommand{\ga}{\gamma}
\newcommand{\si}{\sigma}
\newcommand{\De}{\Delta}
\newcommand{\Om}{\Omega}

\newtheorem{theorem}{Theorem}[section]
\newtheorem{lem}[theorem]{Lemma}

\newproclaim{defn}[theorem]{Definition}
\newproclaim{Remark}{Remark}
\newproclaim{Remarks}{Remarks}
\newproclaim{NR}{Notational Remark}

\newtheorem{obs}[theorem]{Observation}

\makeatother

\setattribute{abstract}   {width}  {290pt}

\begin{document}
\begin{frontmatter}

\title{Sublinearity of the travel-time variance for~dependent first-passage percolation}
\runtitle{Sublinearity for dependent first-passage percolation}

\begin{aug}
\author[A]{\fnms{Jacob} \snm{van den Berg}\corref{}\ead[label=e1]{J.van.den.Berg@cwi.nl}}
and
\author[B]{\fnms{Demeter} \snm{Kiss}\thanksref{t2}\ead[label=e2]{D.Kiss@cwi.nl}}
\runauthor{J. van den Berg and D. Kiss}
\affiliation{CWI and VU University Amsterdam, and CWI}
\address[A]{CWI\\
Science Park 123\\
1098 XG Amsterdam\\
The Netherlands\\
and\\
Department of Mathematics\\
VU University---Faculty of Sciences\\
De Boelelaan 1081a\\
1081 HV Amsterdam\\
The Netherlands\\
\printead{e1}} 
\address[B]{CWI\\
Science Park 123\\
1098 XG Amsterdam\\
The Netherlands\\
\printead{e2}}
\end{aug}

\thankstext{t2}{Supported by NWO.}

\received{\smonth{7} \syear{2010}}
\revised{\smonth{10} \syear{2010}}

%
\begin{abstract}
Let $E$ be the set of edges of the $d$-dimensional cubic lattice
$\Z^d$, with $d \geq2$, and let $t(e), e \in E$, be nonnegative
values. The passage time from a vertex $v$ to a vertex $w$ is defined
as $\inf_{\pi\dvtx v \ra w} \sum_{e \in\pi} t(e)$, where the~%
infi\-mum is over all paths $\pi$ from $v$ to $w$, and the sum is over
all edges $e$ of~$\pi$.\looseness=-1

Benjamini, Kalai and Schramm \cite{BeKaSc} proved that if the $t(e)$'s
are i.i.d. two-valued positive random variables, the variance of the
passage time from the vertex $0$ to a vertex $v$ is sublinear in the
distance from $0$ to $v$. This result was extended to a large class of
independent, continuously distributed $t$-variables by Bena\"{i}m and
Rossignol \cite{BeRo}.

We extend the result by Benjamini, Kalai and Schramm in a very
different direction, namely to a large class of models where the
$t(e)$'s are dependent. This class includes, among other interesting
cases, a~model studied by Higuchi and Zhang \cite{HiZh}, where the
passage time corresponds with the minimal number of sign changes in a
subcritical ``Ising landscape.''
\end{abstract}

%
\begin{keyword}[class=AMS]
\kwd[Primary ]{60K35}
\kwd[; secondary ]{82B43}.
\end{keyword}
\begin{keyword}
\kwd{First-passage percolation}
\kwd{influence results}
\kwd{greedy lattice animals}
\kwd{Ising model}.
\end{keyword}

\end{frontmatter}

\section{Introduction and statement of results}
Consider, for $d \geq2$, the $d$-di\-mensional lattice $\Z^d$. Let $\E
$ denote the set of edges of the lattice, and let~$t(e)$, $e \in\E$,
be nonnegative real values. A path from a vertex $v$
to a~vertex~$w$ is an alternating sequence of
vertices and edges
\[
v_0 = v, e_1, v_1, e_2,\ldots, v_{n-1}, e_n, v_n = w,
\]
where each $e_i$ is an edge between the vertices
$v_{i-1}$ and $v_i$, $1 \leq i \leq n$. To indicate that $e$ is an edge
of a path $\pi$, we often write, with some abuse of notation,
$e \in\pi$.

If $v = (v_1,\ldots, v_d)$ is vertex, we use the notation $|v|$ for
$\sum_{i=1}^d |v_i|$. The (graph) distance $d(v,w)$ between vertices
$v$ and $w$ is defined as $|v - w|$. The vertex $(0,\ldots,0)$ will
be denoted by $0$.\vadjust{\goodbreak}

The passage time of a path $\pi$ is defined as
%
%
\begin{equation} \label{def-pt}
T(\pi) = \sum_{e \in\pi} t(e).
\end{equation}

The passage time (or travel time) $T(v,w)$ from a vertex $v$ to a
vertex $w$ is defined as
\[
T(v,w) = \inf_{\pi\dvtx v \rightarrow w} T(\pi),
\]
where the infimum is over all paths $\pi$ from $v$ to $w$.

Analogous to the above described bond version, there is a natural
site version of these notions: in the site
version, the $t$-variables are assigned to the vertices instead of the
edges. In the definition of $T(\pi)$, the right-hand side in
(\ref{def-pt}) is then replaced by its analog where the sum is over
all vertices of $\pi$. There seems to be no ``fundamental''
difference between the bond and the site version.

An important subject of study in first-passage percolation is the
asymptotic behavior of $T(0,v)$ and it fluctuations,
when $|v|$ is large and the $t(e)$'s are random variables.
It is believed that, for a large class of distributions of the
$t(e)$'s, the variance of $T(0,v)$ is of order
$|v|^{2/3}$. However, this has only been proved for a special case in a
modified (oriented) version of the model~\cite{Jo}.
Apart from this, the best upper bounds obtained for the variance before
2003 were linear in $|v|$ \cite{Ke}. See Section 1 of
\cite{BeKaSc} for more background and references.

Benjamini, Kalai and Schramm \cite{BeKaSc} showed that if the $t(e)$'s
are i.i.d. random variables taking values $a$ and $b$,
$b \geq a > 0$, then the variance of~$T(0,v)$ is sublinear in the
distance from $0$ to $v$.
More precisely, they showed the following theorem.
\begin{theorem}[(Benjamini, Kalai and Schramm
\cite{BeKaSc})]\label{thm-bks}
Let $b \geq a > 0$. If the $(t(e), e \in\E)$ are i.i.d. random
variables taking values in $\{a,b\}$, then
there is a constant $C > 0$ such that, for all $v$ with $|v| \geq2$,
%
%
\begin{equation}\label{ineq-bks}
\operatorname{Var}(T(0,v)) \leq C \frac{|v|}{{\log}|v|}.
\end{equation}
\end{theorem}

Bena\"{i}m and Rossignol \cite{BeRo} extended this result
to a large class of i.i.d. $t$-variables with a continuous
distribution, and also
proved concentration results. See also \cite{Gra}.

We give a generalization of Theorem \ref{thm-bks} in a very different
direction, namely to a large
class of \textit{dependent} $t$-variables. The description of this class,
and the statement of our general results
are given in Section \ref{sect-main-result}.

Using our general results, we show in particular that (\ref{ineq-bks})
holds for
the $\{a, b\}$-valued Ising model with $0 < a < b$ and
inverse temperature $\beta< \beta_c$. By $\{a, b\}$-valued
Ising
model, we mean the model that is simply obtained from
the ordinary,\vadjust{\goodbreak}
$\{-1, +1\}$-valued, Ising model by replacing $-1$ by $a$ and~$+1$ by~$b$.
The precise definition of the Ising model and the statement of this
result is given in Section~\ref{subs-ising}.

We also study, as a particular case of our general results, a different
Ising-like first-passage
percolation model: consider
an ``ordinary'' Ising model (with signs $-1$ and $+1$), with parameters
$\beta< \beta_c$ and with external
field~$h$ satisfying certain conditions. Now define the passage time
$T(v,w)$ between two vertices $v$ and $w$ as
the minimum number of sign changes
needed to travel from $v$ to $w$. Higuchi and Zhang \cite{HiZh}
proved, for $d=2$, a concentration result for this model.
This concentration result implies an upper bound for the variance that
is (a ``logarithmic-like'' factor) \textit{larger}
than linear. We show from our general
framework that the sublinear bound (\ref{ineq-bks}) holds (see Theorem
\ref{thm-ising-sc}).

The last special case we mention explicitly is that where the
collection of $t$-variables is a
finite-valued Markov random field which satisfies
a high-noise condition studied by H\"{a}ggstr\"{o}m and Steif (see
\cite{HaSt}). Again it follows from
our general results that the sublinear bound (\ref{ineq-bks}) holds
(see Theorem \ref{thm-hn}).

The general organization of the paper is as follows: in the next three
subsections, we give precise definitions
and statements concerning the special results mentioned above.
Then, in Section \ref{sect-main-result}, we state our main,
more general results, Theorems \ref{mainthm} and \ref{mainthm-0}.

In Section \ref{proofs-ising-etc}, we prove the special cases
(Theorems \ref{thm-ising},
\ref{thm-hn} and \ref{thm-ising-sc}) from Theorems \ref{mainthm} and
\ref{mainthm-0}.

In Section \ref{section-tools}, we present the main ingredients for
the proofs of our general results:
an inequality by Talagrand (and its extension
to multiple-valued random variables), a very general ``randomization
tool'' of Benjamini, Kalai and Schramm,
and a result on greedy lattice animals by Martin \cite{Ma01}.

In Section \ref{proof-main}, we first give a very brief informal
sketch of the proof of Theorem~\ref{mainthm} (pointing out
the extra problems that arise, compared with the i.i.d. case in \cite
{BeKaSc}), followed by
a formal, detailed proof.

The proof of Theorem \ref{mainthm-0} is very similar to that of
Theorem \ref{mainthm}. This
is explained in Section \ref{proof-mainthm-0}.

\subsection{The case where the $t$-variables have an $\{a, b\}$-valued
Ising distribution}
\label{subs-ising}
Recall that the Ising model (with inverse temperature $\be$ and
external field $h$) on a countably
infinite, locally finite graph $G$ is
defined as follows. First some notation: we write $v \sim w$ to
indicate that two vertices $v$ and~$w$
share an edge. For each vertex $v$ of $G$, the set of vertices $\{v
\dvtx w \sim v\}$ is
denoted by $\partial v$. The spin value ($+1$ or $-1$) at a vertex $v$
is denoted by $\si_v$.
Now define, for each vertex $v$ and each $\al\in\{-1, +1\}^{\partial v}$,
the distribution $q_v^{\al} = q_{v; \be,h}^{\al}$, on $\{-1, +1\}$:
%
%
\begin{eqnarray} \label{def-ising-q}
q_v^{\al}(+1) &=& \frac{\exp(\be(h + \sum_{w \sim v} \al
_w))}
{\exp(\be(h + \sum_{w \sim v} \al_w)) + \exp
(-\be(h + \sum_{w \sim v} \al_w))}, \nonumber\\[-8pt]\\[-8pt]
q_v^{\al}(-1) &=& \frac{\exp(-\be(h + \sum_{w \sim v} \al
_w))}
{\exp(\be(h + \sum_{w \sim v} \al_w)) + \exp
(-\be(h + \sum_{w \sim v} \al_w))}.\nonumber
\end{eqnarray}

Let $V$ denote the set of vertices of $G$.
An Ising distribution on $G$ (with parameters $\be$ and $h$) is a
probability distribution
$\mu_{\be, h}$ on $\{-1, +1\}^V$ which satisfies, for each vertex $v$
and each $\et\in\{-1,+1\}$,
%
%
\begin{equation} \label{def-ising-d}
\mu_{\be, h}(\si_v = \eta| \si_w, w \neq v) = q_v^{\si
_{\partial v}}(\eta),\qquad \mu_{\be, h} \mbox{-a.s.}
\end{equation}

In this (usual) setup, the spin values are assigned to the vertices.
One can define an Ising model with
spins assigned to the edges, by replacing $G$ by its cover graph
(i.e., the graph whose
vertices correspond with the edges of~$G$, and where two vertices share
an edge if the edges of $G$ to
which these vertices correspond, have a common endpoint).

In the case where $G$ is the $d$-dimensional cubic lattice $\Z^d$,
with $d \geq2$, it is well known that there is a critical
value $\be_c \in(0, \infty)$ such that the following holds:
if $\be< \be_c$, there is a unique distribution satisfying (\ref
{def-ising-d}).
If $\be> \be_c$ and $h = 0$, there is more than one distribution
satisfying (\ref{def-ising-d}).
A similar result (but with a different value of $\be_c$) holds for the
edge version of the model.

Let $b > a > 0$. An \textit{$\{a, b\}$-valued Ising model} is obtained
from the usual Ising model by
reading $a$ for $-1$ and $b$ for $+1$. More precisely, if $(\si_v, v
\in V)$ has an Ising distribution
and, for each $v \in V$, $t(v)$ is defined to be $a$ if $\si_v = -1$
and $b$ if $\si_v = +1$,
then we say that $(t(v), v \in V)$ are $\{a, b \}$-valued Ising variables.
A similar definition holds for the situation where the spins are
assigned to the edges.

A special case of our main result is the following extension of Theorem~\ref{thm-bks} to
the Ising model.
\begin{theorem} \label{thm-ising}
Let $b > a >0$ and $d \geq2$.
If $(t(v), v \in\Z^d)$ are $\{a,b\}$-valued Ising variables
with inverse temperature $\be< \be_c$
and external field $h$,
then there is a constant $C > 0$ such that for all $v$ with $|v| \geq2$,
%
%
\begin{equation}\label{ineq-ising}
\operatorname{Var}(T(0,v)) \leq C \frac{|v|}{{\log}|v|}.
\end{equation}
The analog of this result holds for the case where the values $a, b$
are assigned to the edges.
\end{theorem}

\subsection{Markov random fields with high-noise condition}
\label{subs-hn} Let $(\si_v , v \in\Z^d)$, be a translation
invariant Markov random field taking values
in $W^{\mathbb{Z}^d}$ where~$W$ is a finite set.
Let $v \in\Z^d$. For each $w\in W$ define (see \cite{HaSt})
\[
\gamma_w = \min_{\eta\in W^{\partial v}} \mathbb{P}(\si_v =
w | \si_{\partial v} = \eta).
\]

Further, define
\[
\gamma= \sum_{w\in W} \gamma_w.\vadjust{\goodbreak}
\]

Note that the definition of $\gamma_w$ and $\gamma$ does not depend
on the choice of $v$.
H\"{a}ggstr\"{o}m and Steif \cite{HaSt} studied the existence of
finitary codings (and exacts simulations) of
Markov random fields under the
following high-noise (HN) condition (see also \cite{HaNe} and \cite{BeSt}).
\begin{defn}[(HN condition)]\label{HN-def} A translation invariant
Markov random field on $\mathbb{Z}^d$ satisfies the HN condition, if
\[
\gamma>\frac{2d-1}{2d}.
\]
\end{defn}

We will show that the following theorem is a consequence of our main result.
\begin{theorem}\label{thm-hn}
Let $d \geq2$ and let $(\si_v , v \in\Z^d)$ be a translation
invariant Markov random field taking finitely many,
strictly positive values. If this Markov random field satisfies the HN
condition, then, for the
first-passage percolation model with $t(v) = \si_v , v \in\Z^d$,
there is a constant $C > 0$ such that for all $v$ with $|v| \geq2$,
%
%
\begin{equation}\label{ineq-hn}
\operatorname{Var}(T(0,v)) \leq C \frac{|v|}{{\log}|v|}.
\end{equation}
The analog of this result holds for the edge version of the model.
\end{theorem}
\begin{Remark*}
The HN condition for the edge version is a natural modification of that
in Definition \ref{HN-def}. For instance, the $2 d$ in the numerator
and the denominator of the right-hand side of the inequality in
Definition \ref{HN-def} is the number of nearest-neighbor vertices of
a given vertex, and will be replaced by $4 d -2$ (which is the number
of edges sharing an endpoint with a given edge).
\end{Remark*}

\subsection{The minimal number of sign changes in an Ising pattern}
In Section \ref{subs-ising}, the collection of random variables $(t(v),
v \in\Z^d)$ itself had an Ising distribution (with $-1$ and $+1$
translated to $a$, resp., $b$). A quite different first-passage
percolation process related to the Ising model is the one, studied by
Higuchi and Zhang~\cite{HiZh}, where one counts the minimal number of
sign changes from a vertex $v$ to a vertex $w$ in an Ising
configuration.

For $\beta< \beta_c$, let $\theta(\beta, h)$ denote the probability
that $0$ belongs to an infinite $+$ cluster,
and let
\[
h_c(\beta) = \sup\{h \dvtx\theta(\beta,h) = 0\}.
\]
For $d = 2$, it was proved in \cite{Hi2} that $h_c(\beta) > 0$.

Using our general results, we will prove (in Section \ref
{proofs-ising-etc}) the following extension of Theorem
\ref{thm-bks}.
\begin{theorem}\label{thm-ising-sc}
$\!\!$Let the collection of random variables $(\si_v, v \in\Z^2)$ have~an Ising distribution with
parameters $\be< \be_c$ and external field $h$, with
\mbox{$|h| < h_c$}.
Define, for each edge $e= (v_1,v_2)$,
\[
t(e) =
\cases{1, &\quad if $\si_{v_1} \neq\si_{v_2}$,\cr
0, &\quad if $\si_{v_1} = \si_{v_2}$.}
\]
%
For the first-passage percolation model with these $t$-values, there is
a $C > 0$ such that
for all $v$ with $|v| \ge2$,
%
%
\begin{equation}\label{ineq-ising-sc}
\operatorname{Var}(T(0,v)) \leq C \frac{|v|}{{\log}|v|}.
\end{equation}
\end{theorem}
\begin{Remark*}
  Higuchi and Zhang \cite{HiZh} give a
concentration result
for this model (see Theorem 2 in \cite{HiZh}). Their method is very
different from ours.
[It is interesting to note that the paragraph below (1.11) in their
paper suggests that Talagrand-like
inequalities are not applicable to the Ising model.] The upper bound
for the
variance of $T(0,v)$ which follows from their concentration result
is (a ``logarithmic-like'' factor) \textit{larger} than linear. For earlier
results on this
and related models, see the Introduction in \cite{HiZh}.
\end{Remark*}

\subsection{Statement of the main results} \label{sect-main-result}
Our main results, Theorems \ref{mainthm}\break and~\ref{mainthm-0},
involve $t$-variables that can be represented
by (or ``encoded'' in terms of)
i.i.d. finite-valued random variables in a suitable way, satisfying
certain conditions.
These conditions are of the same flavor as (but somewhat different
from) those in Section~2 in \cite{Be08}.

We first need some notation and terminology.
Let $S$ be a finite set, and~$I$ a countably infinite set.
Let $W$ be a finite subset of $I$. If $x \in S^I$, we write~$x_W$ to
denote the tuple $(x_i, i \in W)$.
If $h \dvtx S^I \rightarrow\R$ is a function, and $y \in S^W$, we say
that $y$ determines the
value of $h$ if $h(x) = h(x')$ for all $x$, $x'$ satisfying $x_W = x'_W
= y$.

Let $X_i, i \in I$, be i.i.d. $S$-valued random variables.
We say that the random variables $t(v), v \in\Z^d$, are represented
by the collection $(X_i, i \in I)$,
if, for each $v \in\Z^d$, $t(v)$ is a function of $(X_i, i \in I)$.
The formulation of our main theorems involve certain conditions on such
a representation:

\begin{itemize}

\item
\textit{Condition} (i): There exist $c_0 >0$ and $\ep_0 >0$ such that for
each $v \in\Z^d$
there is a sequence
$i_1(v), i_2(v), \ldots$ of elements of $I$,
such that for all $k = 1, 2,\ldots,$
%
%
\begin{equation}\label{ineq-condi}
P\bigl(\bigl(X_{i_1(v)},\ldots, X_{i_k(v)}\bigr) \mbox{ does not determine }
t(v) \bigr)
\leq\frac{c_0}{k^{3 d + \ep_0}}.
\end{equation}

\item
\textit{Condition} (ii):
%
%
\begin{eqnarray} \label{ineq-condii}
&&\exists\al> 0\ \forall v, w \in\Z^d\ \forall k < \al|v -
w|\nonumber\\[-8pt]\\[-8pt]
&&\qquad\{i_1(v),\ldots, i_k(v)\} \cap\{i_1(w),\ldots, i_k(w)\} =
\varnothing.\nonumber
\end{eqnarray}

\item
\textit{Condition} (iii): The distribution of the family of random
variables $(t(v), v \in\Z^d)$ is
translation-invariant.\vadjust{\goodbreak}
\end{itemize}

We say that the family of random variables $(t(v), v
\in\Z^d)$ has a representation
satisfying conditions (i)--(iii),
if there are $S$, $I$ and i.i.d. $S$-valued random variables $X_i, i
\in I$ as above, such that
the $t$-variables are functions of the $X$-variables satisfying
conditions (i)--(iii) above.

Analogs of these definitions for $t$-variables indexed by the edges of
$\Z^d$ can be given in a straightforward
way.

Now we are ready to state our main theorem.
\begin{theorem} \label{mainthm}
Let $b > a > 0$, and
let, with $d \geq2$, $(t(v), v \in\Z^d)$
be a~family of random
variables that take values in the interval $[a,b]$
and have a~representation satisfying conditions \textup{(i)--(iii)} above.
Then there is a $C > 0$, such that for all $v \in\Z^d$ with $|v| \geq2$,
%
%
\begin{equation}\label{ineq-mainthm}
\operatorname{Var}(T(0,v)) \leq\frac{C |v|}{{\log}|v|}.
\end{equation}
The analog for the bond version of this result also holds.
\end{theorem}

If the $t$-variables can take values equal or arbitrarily close to $0$,
we need a stronger version of condition (i)
and extra condition (iv) (see below).

By an optimal path from $v$ to $w$, we mean a path $\pi$ from $v$ to
$w$ such that $T(\pi) \leq T(\pi')$
for all paths $\pi'$ from $v$ to $w$.

\begin{itemize}

\item
\textit{Condition} (i$'$):
There exist $c_0 >0$, $\ep_0 >0$ and $\ep_1 > 0$, such that for each
$v \in\Z^d$
there is a sequence
$i_1(v), i_2(v), \ldots$ of elements of $I$,
such that for all $k = 1, 2,\ldots,$
%
%
\begin{equation}\label{ineq-condi-prime}
P\bigl(\bigl(X_{i_1(v)},\ldots, X_{i_k(v)}\bigr) \mbox{ does not determine }
t(v) \bigr)
\leq c_0 \exp(- \ep_0 k^{\ep_1}).
\end{equation}

\item
\textit{Condition} (iv):
There exist $c_1, c_2, c_3 >0$ such that for all vertices $v, w$
the probability that there is no optimal path $\pi$ from $v$ to $w$
with $|\pi| \leq c_1 |v - w|$
is at most $c_2 \exp(-c_3 |v - w|)$.
\end{itemize}
\begin{theorem} \label{mainthm-0}
Let $b > 0$, and
let, with $d \geq2$, $(t(v), v \in\Z^d)$
be a collection of random variables
taking values in the interval $[0,b]$,
and having a~representation satisfying conditions \textup{(i$'$), (ii), (iii)} and
\textup{(iv)} above.
Then there is a $C > 0$, such that for all $v \in\Z^d$ with $|v| \geq2$,
%
%
\begin{equation}\label{ineq-mainthm-0}
\operatorname{Var}(T(0,v) \leq\frac{C |v|}{{\log}|v|}.
\end{equation}
The analog of this result for the bond version of the model also holds.
\end{theorem}
\begin{Remarks*}
\begin{longlist}[(b)]
\item[(a)] Note that condition (iii) is in terms of the $t$-variables only: we
do \textit{not} assume that
the index set $I$ has a ``geometric'' structure
and that the $t$-variables are ``computed'' from the $X$-variables in a
``translation-invariant'' way with respect
to that structure (and the structure of $\Z^d$).\vadjust{\goodbreak}

\item[(b)] The goal of our paper is to show that the main result in \cite
{BeKaSc}, although its proof heavily uses
inequalities concerning independent random variables, can be extended
to an interesting class of dependent
first-passage percolation models. In the setup of the above conditions
(i), (ii), (iii), (i$'$) and (iv),
we have aimed to obtain fairly general Theorems \ref{mainthm} and
\ref{mainthm-0}, without becoming
too general (which would give rise to so many extra technicalities that
the main line of argument would be obscured).
For instance, from the proofs it will be clear that there is a kind of
``trade-off'' between conditions (i) and (ii):
one may simultaneously strengthen the first and weaken the second
condition.

Also, if the bound in condition (i$'$) is replaced by a polynomial bound
with sufficiently high degree, Theorem
\ref{mainthm-0} would still hold (but more explanation would be needed
in Section \ref{proof-mainthm-0}).
Since the main motivation for adding this theorem to Theorem \ref
{mainthm} is to handle the interesting
Ising sign-change model studied by Higuchi and Zhang [for which we know
that condition (i$'$) holds] we have not
replaced condition (i$'$) by a weaker condition.
\end{longlist}
\end{Remarks*}

\section{\texorpdfstring{Proofs of Theorems \protect\ref{thm-ising}, \protect\ref{thm-hn} and \protect\ref{thm-ising-sc} from Theorems \protect\ref{mainthm} and \protect\ref{mainthm-0}}
{Proofs of Theorems 1.2, 1.4 and 1.5 from Theorems 1.6 and 1.7}}
\label{proofs-ising-etc}
\subsection{\texorpdfstring{Proof of Theorem \protect\ref{thm-ising}}{Proof of Theorem 1.2}}
In \cite{Be08}, the notion ``nice finitary representation'' has been
introduced in the context of two-dimensional
random fields. See conditions (i)--(iv) in Section 2 of that paper.
In Section 2 (see in particular Theorem~2.3 in that paper), it is shown
that the Ising model with $\be< \be_c$
has such a representation. (See also \cite{BeSt}.) The key ideas and
ingredients are exact simulation by coupling from the past
(see \cite{PrWi} and \cite{BeSt}), and a well-known result by
Martinelli and
Olivieri \cite{MaOl} that under a natural dynamics (single-site
updates; Gibbs sampler) the system has exponential
convergence to the Ising distribution. The random variables used to
execute these updates are taken as the $X$-variables in the
definition of a representation.

Condition (ii) in \cite{Be08} is
somewhat weaker than our current condition (i). However,
as shown in \cite{Be08} (see the arguments between Theorems 2.3 and 2.4
in \cite{Be08}), the above mentioned
exponential convergence shows that the Ising model satisfies an even
stronger bound, namely condition (i$'$) in
our paper.

Condition (iii) in \cite{Be08} corresponds with our condition (ii),
and condition~(iv) in \cite{Be08} is stronger
than our condition (iii).

In \cite{Be08}, only the two-dimensional case is treated (because the
applications are to percolation models where
typical two-dimensional methods are used) but its arguments concerning
``nice finitary representations'' for the Ising
model extend immediately to higher dimensions.

From the above considerations, it follows that the Ising models in the
statement of our Theorem
\ref{thm-ising} indeed have a representation satisfying our conditions
(i)--(iii).
Application of Theorem \ref{mainthm} now gives Theorem \ref{thm-ising}.

\subsection{\texorpdfstring{Proof of Theorem \protect\ref{thm-hn}}{Proof of Theorem 1.4}}
The argument is very similar to that in the proof of Theorem \ref
{thm-ising}. Therefore, we only mention the
points that need extra attention.

As in the proof of Theorem \ref{thm-ising}, the role of the
$X$-variables in Section \ref{sect-main-result}
is played by the i.i.d. random variables driving a single-site update
scheme (Gibbs sampler). In Theorem \ref{thm-ising},
a form of exponential convergence for the Gibbs sampler was used. This
exponential convergence came from
a result in \cite{MaOl}. In the
current situation, the exponential convergence is, as shown in
Proposition 2.1 in \cite{HaSt}, a~consequence of the HN condition. This
exponential convergence implies (again, as in the case of Theorem \ref
{thm-ising}) condition (i) [and, in fact, the
stronger condition (i$'$)] in Section \ref{sect-main-result}. Condition
(iii) is obvious, and condition (ii) follows
easily (as in the proof of Theorem \ref{thm-ising}) from the general
setup of the Gibbs sampler.
So, again, we now apply Theorem \ref{mainthm} to obtain Theorem \ref{thm-hn}.

\subsection{\texorpdfstring{Proof of Theorem \protect\ref{thm-ising-sc}}{Proof of Theorem 1.5}}
Since $\beta< \beta_c$, the collection $(\si_v, v \in\Z^2)$, has
(as pointed out in the proof of Theorem \ref{thm-ising})
a representation satisfying conditions (i), (ii) and (iii). In fact, as
noted in the proof of Theorem \ref{thm-ising},
it even satisfies the stronger form (i$'$) of condition (i). Since
$t(e)$ is a function of the $\si$-values of the two
endpoints of $e$, it follows immediately that the collection $(t(e), e
\in E)$ (where $E$ denotes the set of edges of
the lattice $\Z^2$) satisfies the (bond analog of) the
conditions (i$'$), (ii) and (iii). The fact that (iv) is satisfied
follows immediately from Lemma 6 [and (1.9)]
in \cite{HiZh}.
Theorem \ref{thm-ising-sc} now follows from (the bond version of)
Theorem \ref{mainthm-0}.

\section{\texorpdfstring{Ingredients for the proof of Theorem \protect\ref{mainthm}}{Ingredients for the proof of Theorem 1.6}}
\label{section-tools}
\subsection{An inequality by Talagrand} \label{subs-tal}
Let $S$ be a finite set and $n$ a positive integer. Assign
probabilities $p_s$, $s \in S$, to the elements of $S$.
Let $\mu$ be the corresponding product measure on $\Om:= S^n$.

Let $f$ be a function on $\Om$, and let $\Vert f \Vert_1$ and $\Vert
f \Vert_2$ denote
the $L_1$-norm and $L_2$-norm of $f$ w.r.t.
the measure $\mu$:
\begin{eqnarray*}
\Vert f \Vert_1 &:=& \sum_{x \in\Om} \mu(x) |f(x)|;
\\
\Vert f \Vert_2 &:=& \sqrt{\sum_{x \in\Om} \mu(x)
\vert f (x) \vert^2}.
\end{eqnarray*}

The notation $\bar f_i$ is used for the conditional expectation of $f$
given all
coordinates except the $i$th. More precisely, for $x = (x_1,\ldots,
x_n) \in S^n$ we define
\[
\bar f_i(x) := \sum_{s \in S} p_s f(x_1,\ldots, x_{i-1}, s,
x_{i+1},\ldots, x_n).
\]
Further, we define the function
$\De_i f$ on $\Om$ by
%
%
\begin{equation}\label{de-i-def}
(\De_i f) (x) = f(x) - \bar f_i(x),\qquad x \in\Om.
\end{equation}
\begin{NR*} Often we work with the alternative,
equivalent, description that we have $n$ independent
random variables, say $Z_1,\ldots, Z_n$, with
$P(Z_i = s) = p_s, $ $s \in S, 1 \leq i \leq n$.
To emphasize the identity of the random variables involved, we then
often use the notation
$\De_{Z_i}$ instead of $\De_i$.
\end{NR*}

A key ingredient in \cite{BeKaSc} and in our paper is the following
inequality for the case $|S| = 2$ by Talagrand,
a far-reaching extension of an inequality by Kahn, Kalai and Linial
\cite{KaKaLi}.
\begin{theorem}[(Talagrand \cite{Ta}, Theorem
1.5)] \label{thm-ta}
There is a constant $K > 0$ such that for each $n$ and each function
$f$ on $\{0,1\}^n$,
%
%
\begin{equation} \label{ineq-ta}
\operatorname{Var}(f) \leq K\log\biggl(\frac{2}{p(1-p)}\biggr)
\sum_{i=1}^n \frac{\Vert\Delta_{i}f\Vert_{2}^{2}}
{\log(e\Vert\Delta_{i}f\Vert_{2}/\Vert
\Delta_{i}f\Vert_{1})},
\end{equation}
where (in the notation in the beginning of this section) $p = p_1 = 1 -
p_0$, and
where $\operatorname{Var}(f)$ denotes the variance of $f$ w.r.t. the
measure $\mu$.
\end{theorem}

In the literature, (partial) extensions of this inequality and
inequalities of related flavor, to the case
$|S| > 2$ have been given; see, for example, \cite{Ro08} and~\cite{BeRo}.
The following theorem (see \cite{Ki}) states the most ``literal''
extension of Theorem \ref{thm-ta} to the case $|S| > 2$.
(In \cite{Ki}, an extended version of Beckner's inequality, a key
ingredient in the proof of Theorem \ref{thm-ta}, is used,
and the
proof of Talagrand is followed, with appropriate adaptations, to obtain
the extension of Theorem \ref{thm-ta}.)
To make comparison of our line of arguments with that in \cite{BeKaSc} as
clear as possible, it is this extension we will use. (Moreover, if
instead of Theorem \ref{thm-ki10}
we would use the modified Poincar\'{e}
inequalities in~\cite{BeRo}, this would not simplify our proof of
Theorem \ref{mainthm}.)
\begin{theorem}[(\cite{Ki}, Theorem 1.3)]\label{thm-ki10}
There is a constant $K>0$ such that for each finite set $S$, each $n
\in\N$ and each
function $f$ on $S^n$ the following holds:
%
%
\begin{equation}\label{eq-ki10}
\operatorname{Var}(f) \leq K \log\biggl(\frac{1}{\min_{s \in S}
p_s}\biggr)
\sum_{i=1}^n \frac{\Vert\Delta_{i}f\Vert_{2}^{2}}
{\log(e\Vert\Delta_{i}f\Vert_{2}/\Vert
\Delta_{i}f\Vert_{1})}.
\end{equation}
\end{theorem}

\subsection{Greedy lattice animals} \label{subs-martin}
The subject of this subsection played no role in the treatment of the
first-passage percolation
model with independent $t$-variables in \cite{BeKaSc},
but turns out to be important in our treatment of dependent
$t$-variables.\vadjust{\goodbreak}

Consider, for $d \geq2$, the $d$-dimensional cubic lattice. A lattice
animal (abbreviated as l.a.) is
a finite connected subset of $\Z^d$ containing the origin.
Let~$X_v$, $v \in\Z^d$, be i.i.d. nonnegative random variables with
common distribution~$F$.
Define
\[
N(n) := \max_{\zeta\dvtx\zeta\ \mathrm{l.a.}\ \mathrm{with}\ |\zeta| = n}
\sum_{v \in\zeta} X_v,
\]
where the maximum is over all lattice animals of size $n$.

The subject was introduced by Cox et al. 
\cite{CoGaGrKe}. The asymptotic behavior,
as $n \rightarrow\infty$
of $N(n)$ has been studied in that and several other papers (see, e.g.,
\cite{GaKe} and \cite{HoNe}).
For our purpose, the following result by Martin \cite{Ma01} is very suitable.
\begin{theorem}[(Martin \cite{Ma01}, Theorem
2.3)]\label{thm-martin}
There is a constant $C$ such that for all $n$ and for all $F$ that satisfy
\begin{eqnarray}
&\displaystyle \int_0^{\infty} \bigl(1 - F(x)\bigr)^{1/d} \,d x <
\infty,&
\nonumber\\
%
%
\label{eq-martin}
&\displaystyle E\biggl(\frac{N(n)}{n}\biggr) \leq C \int_0^{\infty}
\bigl(1 - F(x)\bigr)^{1/d} \,d x.&
\end{eqnarray}
\end{theorem}

Martin \cite{Ma01} says considerably more than this, but the above
is sufficient for our
purpose.

\subsection{A randomization tool} \label{randomize}
As in \cite{BeKaSc} we need, for technical reasons, a~certain
``averaging'' argument:
extra randomness is added to the system to make it more tractable. To
handle this extra
randomness appropriately, the following lemma from \cite{BeKaSc} is used.
\begin{lem}[(Benjamini, Kalai and Schramm \cite
{BeKaSc}, Lemma 3)]\label{lem-average}
There is a~constant
$c>0$ such that for every $m\in\mathbb{N}$ there is a~function
\[
g = g_{m}\dvtx\{ 0,1\} ^{m^{2}}\rightarrow\{ 0,1,\ldots
,m\},
\]
which satisfies properties \textup{(i)} and \textup{(ii)} below:

\begin{longlist}
\item For all $i = 1,\ldots, m^2$ and all $x \in\{0,1\}^{m^2}$,
%
%
\begin{equation}\label{eq-g-prop1}
\bigl|g_m\bigl(x^{(i)}\bigr) - g_m(x)\bigr| \leq1,
\end{equation}
where $x^{(i)}$ denotes the element of $\{0,1\}^{m^2}$ that differs
from $x$ only
in the $i$th coordinate.

\item
%
%
\begin{equation}\label{eq-g-prop2}
\max_k \PP\bigl(g(y) = k\bigr) \leq c/m,
\end{equation}
where $y$ is a random variable uniformly distributed on $\{0,1\}^{m^2}$.
\end{longlist}
\end{lem}

\section{\texorpdfstring{Proof of Theorem \protect\ref{mainthm}}
{Proof of Theorem 1.6}} \label{proof-main}
To keep our formulas compact, we will use constants $C_1$, $C_2,
\ldots.$ The precise values of these constants do not
matter for our purposes. Some of them
depend on $a$, $b$, the dimension $d$, the distribution of the
$X$-variables (in terms of which the
$t$-variables are represented), or the constants in the conditions (i),
(i$'$), (ii), (iii) and (iv) in
Section \ref{sect-main-result}. However, they do not (and obviously
should not) depend on the choice of $v$
in the statement of the theorem.

\subsection{Informal sketch} \label{informal}
The detailed proof is given in the next subsection.
Now we first give a very brief and rough summary of the proof of the
main result in \cite{BeKaSc} (listed as
Theorem \ref{thm-bks} in our
paper), and then informally (and again
briefly) indicate the extra problems that arise in our situation where
the $t$-variables are dependent.

Let $\ga$ be the path from $0$ to $v$ for which the sum of the
$t$-variables is minimal. (If more than one such path
exists, choose one of these by a deterministic procedure.)
Since the value of each $t$-variable is at least $a > 0$ and at most
$b$, it is clear that the number of edges of $\ga$ is at most
a constant $c$ times $|v|$.

In \cite{BeKaSc} the $t$-variables are independent, and Talagrand's
inequality (Theorem \ref{thm-ta}) is
applied with $f = T(0,v)$ and with each
$i$ denoting an edge $e$. From the definitions, it is clear that $\De
_i f$ is roughly the change of $T(0,v)$
caused by changing~$t(e)$. Moreover, a change of $t(e)$ can only cause
a change of $T(0,v)$ if, before
or after the change, $e$ is on the above mentioned path $\ga$.
So, ignoring the denominator in Talagrand's inequality, one gets the
(linear) bound
%
%
\begin{equation}\label{bd-lin}
\operatorname{Var}(T(0,v)) \leq C_1 \E\biggl[\sum_{e \in\ga}
(b-a)^2\biggr] \leq c (b-a)^2 |v|.
\end{equation}

It turns out that,
by introducing additional randomness in an appropriate way, without
changing the variance
(see Lemma \ref{lem-average}), the
$\Vert\Delta_{i}f\Vert_{2}/\Vert\Delta
_{i}f\Vert_{1}$ in the denominator in the right-hand side
of Talgrand's inequality becomes (uniformly in $i$) larger than
$|v|^\beta$ for some $\beta> 0$, thus giving
the ${\log}|v|$ (and hence, the sublinearity) in Theorem \ref{thm-bks}.

In our situation, the underlying independent random variables are
the~$X_i$, $i \in I$ (by which the
\textit{dependent} $t$-variables are represented). Application of
Talagrand-type inequalities to these
variables has the complication that changing one $X$-variable changes a
(random) set of possibly many $t$-variables. Taking
the square of the effect complicates this further. Nevertheless, it
turns out that by suitable decompositions
of the summations, and by block arguments (rescaling), one finally
gets, instead of
(\ref{bd-lin})
a bound in terms of (``rescaled'') greedy lattice animals which, by the
result of Martin in
Section~\ref{section-tools}, is still linear in $|v|$.

To handle the denominator in the Talagrand-type inequality, we use
additional randomness, as in \cite{BeKaSc}.
Again, the fact that changing an $X$-variable can have effect on many
$t$-variables complicates the analysis,
but this complication is easier to handle than that for the numerator
mentioned above.

\subsection{Detailed proof} \label{subs-proof-main}
We give the proof for the site version of Theorem~\ref{mainthm}. The
proof for the bond version is obtained from it by
a straightforward, step-by-step translation.
\begin{NR*}
The cardinality of a set $V$ will be indicated by~$|V|$.

We start by stating a simple but important observation (a version of which
was also used in \cite{BeKaSc}). A finite path $\pi$ is called an
optimal path, or a~geodesic, if
there is no path $\pi' \neq\pi$ with the same starting and endpoint
as $\pi$,
for which $T(\pi') < T(\pi)$.
\end{NR*}
\begin{obs} \label{obs-path}
Since the $t$-variables are bounded away from $0$ and~$\infty$,
there is a constant $C_2 > 0$ such that
for every positive integer $n$ and every $w \in\Z^d$ the following
hold:
\begin{longlist}[(a)]
\item[(a)]
Each geodesic has at most $C_2 n$ vertices in the box $w +
[-n,n]^d$.

\item[(b)] Each geodesic which starts at $0$ and ends at $w$ has at most $C_2
|w|$ vertices.
\end{longlist}
\end{obs}
%

Let $X_i, i\in I$, be the independent random variables in terms of
which the
variables $(t(v), v \in\Z^d)$ are represented. So $T(0,v)$ is a
function of the
$X$-variables. As we said in the informal sketch, we introduce extra
randomness, in the same way as in \cite{BeKaSc}: fix $m := \lfloor
|v|^{1/4}\rfloor$.
Let $(y_i^j, i = 1,\ldots, m^2, j = 1,\ldots, d)$ be a
family of independent
random variables, each taking value $0$ or $1$ with probability $1/2$.
The family of $y_i^j$'s is also taken independently of the $X$-variables.
Define, for $j = 1,\ldots, d$,
\[
y^j = (y_1^j,\ldots, y_{m^2}^j).
\]
Each $y^j$ is uniformly distributed on $\{0,1\}^{m^2}$, and will play the
role of the~$y$ in Lemma \ref{lem-average}.
We simply write $Y$ for the collection
$(y_i^j, i = 1,\ldots, m^2, j = 1,\ldots, d)$ and $X$
for the collection $(X_i, i \in I)$.

Let
%
%
\begin{equation}\label{eq-def-zY}
z(Y) = (g(y^1),\ldots, g(y^d))
\end{equation}
with $g = g_m$ as in Lemma \ref{lem-average}.

To shorten notation, we will write $f$ for $T(O,v)$ and $\tilde f$ for
the passage time between the vertices that are obtained from $0$ and $v$
by a (random) shift over the vector $z(Y)$:
%
%
\begin{equation}\label{eq-def-ftil}
\tilde f = T\bigl(z(Y), v + z(Y)\bigr).
\end{equation}

Note that $f$ is completely determined by $X$, while $\tilde f$ depends on
$X$ as well as~$Y$.

By translation invariance [see condition (iii)],
for every $w \in\Z^d$, $T(0,v)$ has the same distribution as $T(w, v
+ w)$.
Hence, by conditioning on $Y$ and using that $Y$ is independent of the
$t$-variables, it follows that $\tilde f$ has the same distribution as $f$.
In particular,
%
%
\begin{equation}\label{compare-f}
\operatorname{Var}(f) = \operatorname{Var}(\tilde f).
\end{equation}

Theorem \ref{thm-ki10} gives (see the Remarks below)
%
%
\begin{eqnarray}\label{appl-tal}
\operatorname{Var}(\tilde{f})
& \leq& C_{3}\sum_{i=1,\ldots,m^{2},j=1,\ldots, d}\frac{
\Vert\Delta_{y_{i}^{j}}\tilde{f}\Vert_{2}^{2}}
{1+\log(\Vert\Delta_{y_{i}^{j}}\tilde{f}\Vert
_{2} {/}
\Vert\Delta_{y_{i}^{j}}\tilde{f}\Vert_{1}
)}\nonumber\\[-8pt]\\[-8pt]
&&{} +C_{3}\frac{\sum_{i \in I}\Vert\Delta_{X_i}\tilde
{f}\Vert_{2}^{2}}
{1+\min_{i \in I}\log(\Vert\Delta_{X_i}\tilde{f}
\Vert_{2} {/}
\Vert\Delta_{X_i}\tilde{f}\Vert_{1})}.\nonumber
\end{eqnarray}
\begin{Remarks*}
\begin{longlist}[(a)]
\item[(a)]
At first sight, Theorem \ref{thm-ki10} is not applicable in the
current situation where
we have two types of random variables: $X_i$'s and $y_i^j$'s. However,
by a straightforward argument,
``pairing'' each variable $y_i^j$, $i = 1,\ldots, m^2$, $j = 1,\ldots, d$,
with an independent ``dummy'' variable $X_i^j$ (with the same
distribution as the ``ordinary'' $X$-variables), and
each variable $X_i$, $i \in I$, with an independent ``dummy'' variable
$y_i$ (with the same distribution as the ``ordinary''
$y$-variables), it is easy to see that Theorem \ref{thm-ki10} is
indeed applicable here.

\item[(b)] Note that the statement of Theorem \ref{thm-ki10} is formulated
for finite $n$. Combined with a standard
limit argument, it gives (\ref{appl-tal}).
\end{longlist}
\end{Remarks*}

We will handle, in separate subsections, the first term of (\ref
{appl-tal}), the numerator
of the second term, and the denominator of the second term.

\subsubsection{\texorpdfstring{The first term in (\protect\ref{appl-tal})}
{The first term in (23)}} \label
{proof-term1}
By (\ref{eq-def-zY}), (\ref{eq-def-ftil}) and (\ref{eq-g-prop1}) it
follows that
$|\De_{y_i^j} \tilde f|$ is at most
a constant $C_{4}$, so that we have the following lemma.
\begin{lem}\label{lem-first-term}
The first term in (\ref{appl-tal}) is at most
%
%
\begin{equation} \label{eq-term1}
\leq d C_{4}m^{2}=d C_{4}|v|^{1/2}.
\end{equation}
\end{lem}

\subsubsection{\texorpdfstring{The denominator of the second term in (\protect\ref{appl-tal})}
{The denominator of the second term in (23)}}
\label{proof-term2-den}
In this subsection we write, for notational convenience, $\Delta_i
\tilde f$ for
$\Delta_{X_i} \tilde f$, where $i \in I$.

If $w, w' \in\Z^d$ we write $\ga_{w,w'}$ for the path $\pi$
minimizing $\sum_{w \in\pi} t(w)$.
If there is more than one such path, we use a deterministic,
translation-invariant
way to select one. If $w = 0$ and $w'$ is our ``fixed'' $v$, we write
simply $\ga$ for $\ga_{0,v}$.\vadjust{\goodbreak}

Recall that $z = z(Y)$ is the random shift.
We write $\ga(z)$ for $\ga_{z, v+z}$.

Also recall the definitions and notation in Section \ref{sect-main-result}.
If $w \in\Z^d$ and \mbox{$j \in I$}, we say that $w$ needs $j$ if $j =
i_k(w)$ for some positive
integer $k$, and $X_{i_1(w)},\ldots, X_{i_{k-1}(w)}$ does not
determine $t(w)$.

By a well-known second-moment argument we have, for each $j \in I$,
%
%
\begin{equation}\label{eq-c-s}
\frac{\Vert\Delta_j \tilde{f}\Vert_{2}}{\Vert
\Delta_j \tilde{f}\Vert_{1}}\geq
\frac{1}{\sqrt{\mathbb{P}(\Delta_j \tilde{f}\neq0 )}}.
\end{equation}

Note that, given $z(Y)$ and all $X_i,i \in I \setminus\{ j \}$,
there is a, possibly nonunique,
$s = s(j,X,Y) \in S$ such that $\tilde{f}$ (now considered as a
function of $X_j$ only)
takes its smallest value at $X_j = s$.
Further note that if $\De_j \tilde f \neq0$ then, after replacing the
value of $X_j$ by $s$, we have
$\De_j \tilde f < 0$. So we get
\[
\mathbb{P}(\Delta_j \tilde{f}<0 )
\geq\mathbb{P}(\Delta_j \tilde{f}\neq0 ) \min
_{r\in S}\mathbb{P}(X_j=r),
\]
and hence
%
%
\begin{equation}\label{eq-25a}
\mathbb{P}(\Delta_j \tilde{f}\neq0 ) \leq
\frac{\mathbb{P}(\Delta_j \tilde{f} < 0 )}
{\min_{r \in S} \mathbb{P}(X_j = r)}.
\end{equation}

Moreover, it follows from the definitions that if $\De_j \tilf< 0$,
there is a $w$ on~$\gamma(z)$ such that
a certain change of $X_j$ causes a change of $t(w)$. By this and~(\ref
{eq-25a}), we have
%
%
\begin{eqnarray} \label{eq-den-1}
\PP(\De_j \tilf\neq0) & \leq& C_5 \sum_{w \in\Z^d} \PP\bigl(w \in
\ga
(z), w \mbox{ needs } j\bigr) \nonumber\\[-8pt]\\[-8pt]
& \leq& C_5 \sum_{w \in\Z^d} \min\bigl(\PP\bigl(w \in\ga
(z)\bigr), \PP(w \mbox{ needs } j)\bigr).\nonumber
\end{eqnarray}

Recall the definition of $m$ in the paragraph following Observation
\ref{obs-path}.
Let $w \in\Z^d$ and consider the box $B_m(w) := w + [-m, m]^d$.
We have
\[
\PP\bigl(w \in\ga(z)\bigr) = \PP(w - z \in\gamma).
\]
By the construction of $z$, and (\ref{eq-g-prop1}), $w-z$ takes values
in the above mentioned
box $B_m(w)$. Also by the construction of $z$, and (\ref{eq-g-prop2}),
each vertex of the box has
probability $\leq C_6/m^d$ to be equal to $w-z$. Moreover, by
Observation~\ref{obs-path} at most $C_7 m$ vertices in the box
are on $\gamma$. Hence, since $\gamma$ is independent of~$z$, it follows
(by conditioning on $\gamma$) that
%
%
\begin{equation}\label{eq-den-2}
\PP\bigl(w \in\gamma(z)\bigr) \leq C_7 m \frac{C_6}{m^d} \leq C_8 |v|^{-(d-1)/4}.
\end{equation}

Further, by condition (i), we have
%
%
\begin{equation} \label{eq-den-3}
\PP(w \mbox{ needs } j) \leq\frac{c_0}{r_w(j)^{3 d + \ep_0}},
\end{equation}
where $r_w(j)$ (which we call the rank of $j$) is the positive integer $k$
for which $i_k(w) = j$.

By (\ref{eq-den-1}), (\ref{eq-den-2}) and (\ref{eq-den-3}), we have,
for every $K$,
%
%
\begin{equation} \label{eq-den-4}
\PP(\De_j \tilf\,{\neq}\,0)\,{\leq}\,C_9 \Biggl(|v|^{-(d-1)/4} |\{w \dvtx
r_w(j) < K\}|\,{+}\,\sum_{k = K}^{\infty} \frac{|\{w \dvtx r_w(j) = k\}|}{k^{3 d +
\ep_0}} \Biggr).\hspace*{-35pt}
\end{equation}

Now, condition (ii) implies, for each $j \in I$ and each $k > 0$,
%
%
\begin{equation}
\label{eq-den-5}
|\{w \dvtx r_w(j) < k\}| \leq C_{10} k^d.
\end{equation}

Hence, the first term between the brackets in (\ref{eq-den-4}) is at most
%
%
\begin{equation}\label{eq-den-6}
C_{10} |v|^{-(d-1)/4} K^d.
\end{equation}

Further, using again (\ref{eq-den-5}) (and summation by parts)
the sum over $k$ in~(\ref{eq-den-4}) is at most
%
%
\begin{equation}\label{eq-den-7}
C_{11} \sum_{k = K}^{\infty} \frac{k^{d}}{k^{3 d + \ep_0 +1}}
\leq
C_{12} K^{-2 d - \ep_0}.
\end{equation}

Combining (\ref{eq-den-4}), (\ref{eq-den-6}) and (\ref{eq-den-7}),
we get
%
%
\begin{equation}\label{eq-den-10}
\PP(\De_j \tilf\neq0) \leq C_{13} \bigl(|v|^{-(d-1)/4} K^d + K^{-2
d - \ep_0} \bigr).
\end{equation}

Now take for $K$ the smallest positive integer satisfying $K^d \geq
|v|^{(d-1)/8}$ and insert this in
(\ref{eq-den-10}). This gives
%
%
\begin{equation} \label{eq-den-11}
\PP(\De_j \tilf\neq0) \leq C_{14} |v|^{-(d-1)/8},
\end{equation}
which together with (\ref{eq-c-s}) yields the following lemma.
\begin{lem} \label{lem-den}
There is a constant $C_{15} > 0$ such that for all $v \in\Z^d$
the denominator of the second term in (\ref{appl-tal}) is larger than
or equal to
\[
{C_{15} \log}|v|.
\]
\end{lem}

\subsubsection{\texorpdfstring{The numerator of the second term in (\protect\ref{appl-tal}), and completion of the proof of Theorem \protect\ref{mainthm}}
{The numerator of the second term in (23), and completion of the proof of Theorem 1.6}}
\label{proof-term2-num}
As in the previous subsection, we write $\De_j$ for $\De_{X_j}$,
where $j \in I$.

By the definition of $\tilf$ (and of the norm $\Vert\cdot\Vert_2$),
we rewrite
%
%
\begin{equation} \label{eq-num-1}
\sum_{j \in I} \Vert\De_j \tilf\Vert_2^2 = \sum_{j \in I} \E
\bigl[\bigl( \De_j T \bigl( z(Y),
z( Y) + v\bigr)\bigr)^2\bigr].
\end{equation}

By taking the expectation outside the summation, conditioning on $Y$
(and using that $Y$ is independent
of the $t$-variables) and then taking the expectation back inside
the\vadjust{\goodbreak}
summation, it is clear that the right-hand side of
(\ref{eq-num-1}) is smaller than or equal to
%
%
\begin{equation} \label{eq-num-2}
\max_{x \in\Z^d} \sum_{j \in I} \E\bigl(\bigl(\De_j T(x,
x+v)\bigr)^2\bigr).
\end{equation}

We will give an upper bound for the sum in (\ref{eq-num-2}) for the
case $x=0$. From the computations, it will
be clear that this upper bound does not use the specific choice of $x$,
and hence holds for all $x$.

In the case $x=0$, the sum in (\ref{eq-num-2}) is, by definition, of course
%
%
\begin{equation}\label{eq-num-2-a}
\sum_{j \in I} \Vert\De_j f \Vert_2^2.
\end{equation}

Let $X'_j$ be an auxiliary random variable that is independent of the
$X$-variables and has the same distribution.
Let $X$ denote the collection of random variables $(X_i, i \in I)$, and
$X'$ the collection obtained from the
collection $X$ by replacing $X_j$ by $X'_j$.
By the definition of $\De_j f$ (and standard arguments), we have
%
%
\begin{eqnarray} \label{gen-var-eq}
\E_j((\De_j f)^2) & = & \tfrac{1}{2}
\E_{j, j'}\bigl[\bigl(f(X) - f(X')
\bigr)^2\bigr] \nonumber\\[-8pt]\\[-8pt]
& = & \E_{j, j'}\bigl[\bigl(f(X) - f (X'
)\bigr)^2
I\bigl(f(X) < f(X')\bigr)\bigr],\nonumber
\end{eqnarray}
where $\E_j$ denotes the expectation with respect to $X_j$, and $\E
_{j,j'}$ denotes the expectation with respect to $X_j$ and
$X'_j$. [So, (\ref{gen-var-eq}) is a function of the collection $(X_i,
i \in I, i \neq j)$.]

Let $\ga$ be the optimal path, as defined in the beginning of
Section \ref{proof-term2-den}, w.r.t. the $t$-variables
corresponding with the family $X$. Let $w$ be a vertex. Observe that a
change of $t(w)$ does not increase $f$ if
$w$ is not on $\ga$, and increases $f$ by at most $b-a$ if $w$ is on
$\ga$. By this observation, and a similar argument as used for
(\ref{eq-den-1}), we have
%
%
\begin{equation} \label{eq-num-4}
\bigl(f(X') - f(X)\bigr) I\bigl(f
(X) < f(X')\bigr) \leq
(b-a) \sum_{w \in\ga} I(w \mbox{ needs } j),
\end{equation}
and hence
%
%
\begin{eqnarray} \label{eq-num-4-a}
&&\bigl(f(X) - f(X')\bigr)^2 I\bigl(f(X) < f(X')\bigr)\nonumber\\[-8pt]\\[-8pt]
&&\qquad \leq(b-a)^2 \sum_{u, w \in\ga}
I(u \mbox{ needs } j, w \mbox{ needs } j).\nonumber
\end{eqnarray}

Since $\Vert\De_j f \Vert_2^2$ is the expectation w.r.t. the $X_i, i
\neq j$, of $\E_j((\De_j f)^2)$, we have, by
(\ref{gen-var-eq}) and (\ref{eq-num-4-a}), that
%
%
\begin{equation}\label{eq-num-4-b}
\Vert\De_j f \Vert_2^2 \leq(b-a)^2 \E\biggl[ \sum_{u, w \in
\ga} I(u \mbox{ needs } j, w \mbox{ needs } j)\biggr].
\end{equation}

To bound the right-hand side of~(\ref{eq-num-4-b}), recall the definition
[below (\ref{eq-den-3})]
of~$r_w(j)$ (with $j \in I$\vadjust{\goodbreak} and $w \in\Z^d$), and
note that, by condition (i) in Section~\ref{sect-main-result},
we have, on an event of probability $1$,
%
%
\begin{eqnarray}\label{eq-num-5}
&& \sum_{u, w \in\ga} I(u \mbox{ and } w \mbox{ need } j)
\nonumber\\[-3pt]
&&\qquad = \sum_{k = 1}^{\infty} \sum_{u, w \in\ga} I\bigl( u \mbox{
and } w \mbox{ need } j,
\max(r_u(j), r_{w}(j)) = k \bigr) \nonumber\\[-9.5pt]\\[-9.5pt]
&&\qquad \leq 2 \sum_{k = 1}^{\infty} \sum_{u \in\ga} \sum_{w \in\ga
} I\bigl( u \mbox{ and } w \mbox{ need } j,
r_u(j) = k, r_{w}(j) \leq k \bigr) \nonumber\\[-3pt]
&&\qquad \leq 2 \sum_{k = 1}^{\infty} \sum_{u \in\ga} I\bigl(u \mbox{ needs
} j, r_u(j) = k\bigr)
|\{w \in\ga\dvtx r_{w}(j) \leq k\}|.\nonumber\vspace*{-2pt}
\end{eqnarray}

By condition (ii), each of the vertices $w$ in the last line of (\ref
{eq-num-5}) is located in a hypercube of length
$C_{16} k$ centered at $u$. By this and Observation \ref{obs-path},
it follows that the number of $w$'s
in the last line of (\ref{eq-num-5}) is at most $C_{17} k$.
So we have, with $C_{18} = 2 C_{17}$,
\[
\sum_{u, w \in\ga} I(u \mbox{ and } w \mbox{ need } j) \leq
C_{18} \sum_{k = 1}^{\infty} k \sum_{u \in\ga} I\bigl(u \mbox{ needs }
j, r_u(j) = k\bigr),\vspace*{-2pt}
\]
which, together with (\ref{eq-num-4-b}) [and using the definition of
$i_k(u)$] gives, after summing over $j$,
%
%
\begin{eqnarray} \label{eq-num-6}
\sum_{j \in I} \Vert\De_j f \Vert_2^2 & \leq& C_{19} \sum_{k =
1}^{\infty} k
\E\biggl[\sum_{u \in\ga} I\bigl(u \mbox{ needs } i_k(u)\bigr)\biggr] \nonumber\\[-3pt]
&=& C_{19} \sum_{k = 1}^{|v|} k \E \biggl[\sum_{u \in\ga} I\bigl(u
\mbox{ needs } i_k(u)\bigr)\biggr]\\[-3pt]
&&{} + C_{19} \sum_{k > |v|} k \E
\biggl[\sum_{u \in\ga} I\bigl(u \mbox{ needs } i_k(u)\bigr)\biggr].\nonumber\vspace*{-2pt}
\end{eqnarray}

The sum over $k > |v|$ in the right-hand side of (\ref{eq-num-6}) can be
bounded very easily
as follows:
by Observation \ref{obs-path}(b), all vertices of $\ga$ are inside
the box
$[- C_2 |v|, C_2 |v|]^d$. Hence, the above-mentioned sum over $k > |v|$
is at most
\[
C_{19} \sum_{k > |v|} k \sum_{u \in[- C_2 |v|, C_2 |v|]^d} \PP\bigl(u
\mbox{ needs } i_k(u)\bigr).\vspace*{-2pt}
\]
By condition (i), and since the number of vertices $u$ in this last expression
is, of course, of order $|v|^d$, this expression is smaller than or equal
to a constant times
\[
|v|^d \sum_{k > |v|} k^{1 - 3 d - \ep_0},\vspace*{-2pt}
\]
which is smaller than a constant\vadjust{\goodbreak} $C_{20}$. 

To bound the sum over $k \leq|v|$ in the right-hand side of (\ref{eq-num-6}),
observe that, by condition (ii), if a set $V \subset\Z^d$ is such
that $|u - u'| \geq C_{21} k$ for all $u, u' \in V$
with $u \neq u'$, then the collection of random variables
\[
\bigl( I\bigl(u \mbox{ needs } i_k(u)\bigr), u
\in V\bigr)
\]
is independent.
With this in mind, we partition, for each $k$, $\Z^d$ in boxes
\[
B_k(w) := [-\lceil C_{21} k \rceil, \lceil C_{21} k \rceil
)^d + 2 \lceil C_{21} k \rceil w,\qquad w \in\Z^d.
\]
We will say that two boxes $B_k(w)$ and $B_k(u)$ are neighbors [where
$u = (u_1,\ldots, u_d)$ and
$w = (w_1,\ldots, w_d)$] if $\max_{1 \leq i \leq d} |w_i - u_i| = 1$.

By Observation \ref{obs-path}(a), $\ga$ has at most $C_{22} k$
vertices in each of these boxes.
Hence, the the sum over $k \leq|v|$ in the right-hand side of (\ref{eq-num-6})
is at most
%
%
\begin{equation} \label{eq-num-7}
C_{23} \sum_{k = 1}^{|v|} k^2
\E\biggl[\sum_{ w \dvtx(*)} I\bigl(\exists u \in B_k(w) \mbox
{ s.t. } u \mbox{ needs } i_k(u)\bigr)\biggr],
\end{equation}
where $(*)$ indicates that we sum over all $w \in\Z^d$ with the
property that $\ga$ has a vertex in $B_k(w)$ or
in a neighbor of $B_k(w)$.

Next, partition $\Z^d$ in $2^d$ classes, as follows:
\[
\Z_z := z + 2 \Z^d,\qquad z \in\{0,1\}^d.
\]
So (\ref{eq-num-7}) can be written as
%
%
\begin{equation} \label{eq-num-8}\qquad
C_{23} \sum_{k = 1}^{|v|} k^2 \sum_{z \in\{0,1\}^d } \E\biggl[\sum
_{w \dvtx(**)}
I\bigl(\exists u \in B_k(z + 2 w) \mbox{ s.t. } u \mbox{ needs }
i_k(u)\bigr)\biggr],
\end{equation}
where $(**)$ indicates that we sum over all $w \in\Z^d$ with the
property that $\ga$ has a point in $B_k(z + 2 w)$ or
in a neighbor of $B_k(z + 2 w)$.

Now, for each $z \in\{0,1\}^d$, the set
\[
\{w \in\Z^d \dvtx\ga\mbox{ has a point in } B_k(z + 2 w) \mbox{
or a neighbor of } B_k(z + 2 w)\}
\]
is a lattice animal and has, for $k\!\leq\!|v|$, by Observation \ref
{obs-path}(b), at most~$C_{24} |v| / k$ elements.

So, from (\ref{eq-num-6})--(\ref{eq-num-8}) we get
%
%
\begin{eqnarray} \label{eq-num-9}
\hspace*{23pt}\sum_{j \in I} \Vert\De_j f \Vert_2^2
&\leq& C_{23} \sum_{k = 1}^{|v|} k^2 \sum_{z \in\{0,1\}^d} \E
\biggl[ \max_{\calL\dvtx| \calL| \leq C_{24} |v|/k} \sum_{w \in
\calL} I\bigl(\exists u \in B_k(z + 2 w) \nonumber\\
\hspace*{23pt}&&\hspace*{179.2pt}\mbox{ s.t. } u \mbox
{ needs } i_k(u)\bigr)\biggr] \\
\hspace*{23pt}&&{} + C_{20},\nonumber
\end{eqnarray}
where the maximum is over all lattice animals $\calL$ with size\vadjust{\goodbreak}
$\leq
C_{24} |v|/k$.

Now for each $z$ we have, by the observation below (\ref{eq-num-6}), that
\[
\bigl(I\bigl(\exists u \in B_k(z + 2 w) \mbox{ s.t. } u
\mbox{ needs } i_k( u)\bigr), w \in\Z^d \bigr)
\]
is a collection of independent 0--1-valued random variables. For each
$w$, this random variable is
$1$ with probability less than or equal to
%
%
\begin{equation}\label{eq-num-10}
\vert B_k(z + 2 w)\vert\max_{u \in\Z^d} \PP\bigl(u \mbox{ needs }
i_k(u)\bigr) \leq\frac{C_{25} k^d}{k^{3 d + \ep_0}},
\end{equation}
where we used condition (i).

By (\ref{eq-num-9}), (\ref{eq-num-10}) and Theorem
\ref{thm-martin}, we get
%
%
\begin{eqnarray}\label{eq-num-11}
\sum_{j \in I} \Vert\De_j f \Vert_2^2 & \leq& C_{20} +
C_{26} \sum_{k = 1}^{|v|} k^2 \frac{|v|}{k} \biggl(\frac{k^d}{k^{3 d
+ \ep_0}}\biggr)^{{1}/{d}} \nonumber\\
& \leq& C_{20} + C_{26} |v| \sum_{k = 1}^{\infty} k^2 k^{-(3 d + \ep
_0)/d} \\
& \leq& C_{27} |v|.\nonumber
\end{eqnarray}

Together with (\ref{eq-num-1})--(\ref{eq-num-2-a}), this gives the
following lemma.
\begin{lem}\label{lem-num-sec}
The numerator of the second term in (\ref{appl-tal}) is at most
$C_{27} |v|$.
\end{lem}

Lemma \ref{lem-num-sec}, together with (\ref{compare-f}), (\ref
{appl-tal}), Lemma \ref{lem-first-term} and
Lemma \ref{lem-den}, completes the proof of Theorem \ref{mainthm}.

\section{\texorpdfstring{Proof of Theorem \protect\ref{mainthm-0}}
{Proof of Theorem 1.7}} \label{proof-mainthm-0}
The proof is very similar to that of Theorem~\ref{mainthm} and we only
discuss those steps that need
adaptation.

First, we define, for $u, w \in\Z^d$,
the following modification of $T(u,w)$:
%
%
\begin{equation} \label{eq-Th}
\hT(u,w) := \min_{\pi\dvtx u \rightarrow w, |\pi| \leq c_1 |u - w|}
T(\pi),
\end{equation}
where $|\pi|$ is the number of vertices of $\pi$, and with $c_1$ as
in condition (iv).

From this definition, it is obvious that
\begin{eqnarray*}
&&|\hT(0,v) - T(0,v)| \\
&&\qquad\leq b (|v| +1) I(\nexists\mbox{ an
optimal path }
\pi\mbox{ from } 0 \mbox{ to } v \mbox{ with } |\pi| < c_1 |v|).
\end{eqnarray*}

By this inequality and condition (iv), we get immediately
\[
\operatorname{Var}(\hT(v)) - \operatorname{Var}(T(v)) = o\bigl(|v|/\log(|v|)\bigr),
\]
so that it is sufficient to prove (\ref{ineq-mainthm-0}) for $\hT(0,v)$.

Now, with $f = \hT(0,v)$ and $\tilde f = \hT(z, v+z)$ [with $z =
z(Y)$ as
in Section \ref{proof-main}] the proof follows that of Theorem \ref
{mainthm}, with the
following modifications:

A few lines above (\ref{eq-den-2}) we used that
$\ga$ has at most $C_7 m$ vertices in the box~$B_m(w)$.
In the current situation we have to add, as a correction term, the probability
that $\ga$ has more than $C_7 m$ vertices in that box. It follows
easily from condition (iv)
that, with a proper choice of $C_7$, this probability goes to $0$
faster than any power of $m$.
Hence (recalling the definition of $m$), it is clear that (\ref
{eq-den-2}) remains true.
Therefore, the denominator of the second term in the proof of Theorem \ref{mainthm}
is, in the current situation, again larger than a constant times ${\log}|v|$.

A few lines before (\ref{eq-num-6}) we applied Observation \ref
{obs-path}(a) (which used
the fact that all $t$-values were larger than some positive $a$) to
conclude that
the number of vertices of $\ga$ in a certain box of length of order
$k$ is at most some constant
times $k$. In the current situation we do not have this strong bound,
but we can obviously conclude that
this number is at most the total number of vertices in the box.
Because of this, the $k$ in (\ref{eq-num-6}) is, in our current
situation, replaced
by $k^d$.

A few lines above (\ref{eq-num-7}), we again used Observation \ref
{obs-path}(a).
Again we have to replace a factor $k$ by $k^d$. By this (and the
previous remark) the $k^2$ in~(\ref{eq-num-7}), and therefore also in (\ref{eq-num-8}) becomes
$k^{2 d}$.

By the definition of $\hT$, the statement about the size of the
lattice animal [a few lines above
(\ref{eq-num-9})] still holds (with appropriate constants).
By this and the earlier remarks, we now get (\ref{eq-num-9}) with the
factor $k^2$ replaced by $k^{2 d}$.
By condition (i$'$), the denominator in the right-hand side of (\ref{eq-num-10})
is now of order $\exp(\ep_0 k^{\ep_1})$,
so that the sum over $k$ in this modified form of (\ref{eq-num-11}) is
still finite.

This completes the proof of Theorem \ref{mainthm-0}.

\section*{Acknowledgments}
Our interest in this subject was triggered by a lecture series on
first-passage percolation by Vladas Sidoravicius in the fall of 2009.

We thank Chuck Newman for referring us to certain results on greedy
lattice animals.

%

%
\printaddresses

\end{document}